\documentclass[11pt]{amsart} 

 \usepackage{pdfsync}
\usepackage{lmodern}
\usepackage{amsthm,amsmath,amssymb,amscd,graphicx,enumerate, stmaryrd,xspace,verbatim, epic, eepic,color,url}
\usepackage{tikz}
  \usepackage{pgfplots}
\pgfplotsset{compat=1.16}

 \usepackage{tikz-3dplot}
\usepackage{hyperref}
\usepackage{pgf}

   \newtheorem{thm}{Theorem}[subsection]
      \newtheorem*{thm*}{Theorem}

   \newtheorem{prop}[thm] {Proposition}     
   \newtheorem{lemma} [thm]{Lemma}
   \newtheorem*{claim}{Claim}

   \newtheorem{cor} [thm]{Corollary}
   \newtheorem{conjecture}{Conjecture}

   \newtheorem*{conjecture*}{Conjecture}
   \newtheorem*{remark*}{Conjecture}

\theoremstyle{definition}
 
          \newtheorem*{exercise*}{Exercise}

  \newtheorem{defi}[thm] {Definition}

  \newtheorem{remark} [thm]{Remark}

\usepackage{amsmath}
\usepackage{amssymb,amsthm}
\usepackage{amscd}
\usepackage{enumerate}
\input{xy}
\xyoption{all}
  
\newcounter{prg}[section]

\newcommand{\la}{\longrightarrow}

\newcommand{\C}{\mathbb{C}}

\newcommand{\N}{\mathbb{N}}
\newcommand{\PP}{\mathbb{P}}

\newcommand{\cE}{\mathcal {E}}

\newcommand{\Hyp}{\operatorname{Hyp}}

\newcommand{\Pic}{\operatorname{Pic}}

\newcommand{\mult}{{\operatorname{mult}}}

\newcommand{\n}{n}
\begin{document}
 \title[Hypertangency and the exceptional set]{Hypertangency of plane curves and  the algebraic exceptional set} 
  \author{Lucia Caporaso and Amos Turchet}
  \address[Caporaso]{Dipartimento di Matematica e Fisica\\ Universit\`{a} Roma Tre \\ Largo San Leonardo Murialdo \\I-00146 Roma\\  Italy }\email{lucia.caporaso@uniroma3.it}
  \address[Turchet]{Dipartimento di Matematica e Fisica\\ Universit\`{a} Roma Tre \\ Largo San Leonardo Murialdo \\I-00146 Roma\\  Italy }\email{amos.turchet@uniroma3.it}
\begin{abstract}  We  investigate plane curves intersecting in at most two unibranch points  to study the algebraic exceptional set appearing in standard conjectures of diophantine and hyperbolic geometry. Our  first result compares the local geometry of two hypertangent curves, i.e. curves having maximal contact at one unibranch point. This is applied to fully describe the exceptional set and, more generally, the hyper-bitangency set, of a plane curve with three components. \end{abstract}
\maketitle

\tableofcontents
 
\section{Introduction}
The goal of this paper is to study  complex plane curves that intersect a given curve in few points.
The main motivation is the quasi-projective analogue of a well-known conjecture of Lang  predicting that on a projective surface of general type there are only finitely many rational and elliptic curves. The  analogue for the complement of  a curve in $\PP^2$ is the following.

\begin{conjecture}\label{conj:main}
  Let $B \subset \PP^2$ be a reduced curve with at most   normal crossing singularities (i.e. nodes). If $\deg B \geq 4$, then the set of   rational curves $C\subset \PP^2$ such that $\vert \nu_C^{-1} (C \cap B) \vert \leq 2$ is finite, where $\nu_C: C^\nu \to C$ is the normalization.
\end{conjecture}

The  set in the conjecture is often referred to as the \emph{algebraic exceptional set}  (see for example \cite{Lan86}); we denote it by $\cE(B)$, i.e. 
\begin{equation}
 \label{eq-EB}
\cE(B):=\{C\subset \PP^2:\ C \text{ rational}, \  \  |\nu_C^{-1}(C\cap B)|\leq 2\}. 
\end{equation}

The curves $C \setminus B$ with $C \in \cE(B)$ are quasi-projective analogues of rational and elliptic curves, both arithmetically (they have a potentially infinite set of integral points when defined over a number field) and geometrically (see for example \cite{KMcK}).

Conjecture~\ref{conj:main} can be deduced by several conjectures of Lang, Vojta, Demailly and Campana, as in \cite[Conj. 3.7]{Lang91}, \cite[Conj. 3.4.3]{Vojta87} \cite{De97} and \cite[Conj. 9.2 and 9.20]{Ca04} (see \cite{BG,SH,AT_book,DT_invit} for introductions and discussions), on the distribution of integral and rational points on quasi-projective surfaces defined over number fields and function fields. From the arithmetic point of view, $\cE(B)$ is conjecturally made of the curves that contain all  potentially infinite families of integral points in the affine surface $\PP^2 \setminus B$.

Historically, describing the set $\cE(B)$, has proven harder the fewer the irreducible components of the curve $B$, and Conjecture~\ref{conj:main} is open if $B$ has at most two components. When $B$ has at least three  components, Corvaja and Zannier proved the function field version of Vojta's Conjecture for  $\PP^2 \setminus B$,  in \cite{CZConic} and \cite{CZGm}, which implies that  curves in $\cE(B)$ have bounded degree (the same result in the so-called non-split case was proved in \cite{CaTur}, building upon \cite{Tur}). More recently \cite{Wang_etal} showed (as a special case of their Theorem 4) that the  set $\cE(B)$ corresponds to  a closed subset of $\PP^2$ and described a closed set containing it. If two of the three components are lines, this can be deduced from an earlier result in \cite{CZ00}, which covers also some non normal-crossing $B$, see Subsection \ref{sec-snc}  for further details.

A different line of investigation assumes that   the curve $B$ is very general. In  \cite{CRY} the authors showed that when $B$  is very general of degree 4 the set $\cE(B)$ consists only of the bitangents and the flex lines, answering a question of Lang. On the other hand, when $B$  is very general of degree at least 5 the  set $\cE(B)$ is empty: this has been proven independently in \cite{Chen} and in \cite{PR}. No similar result is known without the ``very general" assumption.  

In the present article  we extend our analysis from rational curves to curves of arbitrary genus.
This enables us, in particular, to
 give an explicit description of the set $\cE(B)$ when
 $B$ has (at least) three irreducible components. 
 To state our   results   we introduce, generalizing  \eqref{eq-EB},  the set of curves  ``hyper-bitangent" to $B$ 
$$
\Hyp(B,2):=\{C\subset \PP^2:\ C \text{ integral}, \  \  |\nu_C^{-1}(C\cap B)|\leq 2\}  
$$
so that, of course, $\cE(B)\subset \Hyp(B,2)$.
We now summarize in one statement our main results; see Proposition~\ref{prop-3C1}, and Theorems~\ref{3C},   \ref{main}, and  \ref{thm-3Cgen} for more precise and stronger statements.
\begin{thm*}
    Let $B \subset \PP^2$ be a  reduced curve with three  irreducible components and  only nodal   singularities.
If $\deg B\geq 4$ then
\begin{enumerate}[(a)]
 \item
 $\cE(B)= \Hyp(B,2)$;
\item
$ \Hyp(B,2)$ is finite and effectively bounded;
\item
if $\deg B\geq 5$ and $B$ is general, then 
  $\Hyp(B,2)$ is empty.
  \end{enumerate}
  \end{thm*} 
  
  The above effective  bounds on $|\Hyp(B,2)|$ depend only on the   degrees of the irreducible components of $B$. 
The phrase   ``$B$ is general"    means that, for every fixed triple of degrees summing to $d$, the curve $B$ varies in an open dense subset of the space of  triples of  curves of  fixed degrees.
  
  Proposition~\ref{prop-3C1} gives the bound for the number of  hyper-bitangent lines, and Theorem~\ref{thm-3Cgen} gives the bound for  the number of  hyper-bitangent curves of degree at least $2$.

The  proof of these facts relies on a geometric result of independent interest  proved at the beginning of the paper, Theorem~\ref{mirror-sm}, which analyzes two  plane curves, $B$ and $C$, meeting in only one point. 
 The precise statement requires some technical preliminaries;  informally speaking, it establishes that the local geometries of $B$ and $C$ at their intersection point are closely related.

Our results are related to the hyperbolic properties of $\PP^2 \setminus B$. Demailly in \cite{De97} introduced an algebraic analogue of  hyperbolicity for projective varieties, which was extended to quasi projective varieties  in \cite{Chen}. Using this analogue, if  $B$ has degree at least 4   one expects the existence of a positive constant $A$ such that, for every integral curve $C \subset \PP^2$ not contained in the exceptional set $\cE(B)$, the following  bound holds:
\begin{equation}\label{eq:alg-hyp}
  \deg C \leq A \cdot \left( 2g(C^\nu) - 2 + \vert \nu^{-1}_C(C \cap B) \vert \right).
\end{equation}
Equivalently one expects that $\PP^2 \setminus B$ is ``algebraically hyperbolic'' modulo $\cE(B)$; see \cite[Section 9]{J_book} for further  references and discussions on algebraic hyperbolicity.

To see the link with our results,  curves in $\Hyp(B,2)$ correspond exactly to curves that satisfy $\vert \nu_C^{-1}(C \cap B)\vert \leq 2$.
Hence our Theorem proves that, when $B$ has at least three irreducible components and $\vert \nu_C^{-1}(C \cap B)\vert = 2$, the bound \eqref{eq:alg-hyp} holds. More precisely, we show that the degree of a curve in $\Hyp(B,2)$ is  bounded uniformly independently of the genus, and  \eqref{eq:alg-hyp} holds for curves in $\Hyp(B,2)$ with $A = \deg B - 2$; see Corollary \ref{cor:Vojta_bound}. This strengthens, in this particular case,  \cite[Theorem 1]{CZGm} which proves (weak) algebraic hyperbolicity of $\PP^2 \setminus B$ when $B$ has at least three components. We point out that our method is completely different and  uses purely geometric techniques. 

In a parallel direction, the logarithmic Green-Griffiths-Lang Conjecture predicts that, when $\deg B \geq 4$, there exists an \emph{analytic exceptional set}, i.e. a Zariski proper
 closed subset of $\PP^2 \setminus B$ that contains all the images of the non constant holomorphic maps $\C \to \PP^2 \setminus B$, i.e. $\PP^2 \setminus B$ is Brody hyperbolic modulo such analytic exceptional set.
In general, the analytic exceptional set always contains the algebraic one, and they are conjectured to coincide, by Lang.
 In our setting the two sets are known to coincide, and hence our results provide the description of the analytic exceptional set for $\PP^2 \setminus B$, in the case where $B$ has at least three irreducible components. This is only possible since algebraic degeneracy (i.e. every holomorphic map $\C \to \PP^2 \setminus B$ has non dense image) was already proven in \cite{NWY}.

   Our arguments do not use  any of the known results  we mentioned in this introduction.
Therefore
we expect   our techniques to extend to  log  surfaces  where   boundedness is not known.

 \subsection*{Outline of the paper.}
 In Section~\ref{sec-hyp} we establish our principal geometric tool, Theorem~\ref{mirror-sm}, on the local geometry of two plane curves meeting in only one  unibranch point.
 In  Section~\ref{sec-3C} we prove our main results describing the sets $\Hyp(B,2)$ for curves with three components, thus proving the previously stated Theorem. Finally, Section~\ref{sec-examples} collects various special cases and examples related to the earlier topics.

\subsection*{Notation and terminology} We work over $\C$.
We  denote by $C$
   an integral (i.e. reduced and irreducible) projective curve lying in a smooth projective   surface $S$,
  by $\nu_C:C^{\nu}\to C$ its normalization, 
by $p_a(C)$  its arithmetic genus, and  by  $g(C)=p_a(C^{\nu})$ its    geometric genus.
A {\it rational} curve is an integral curve  of geometric genus zero.
Given a point   $p\in C$,
we write $\mult_p(C)$ for the multiplicity of $C$ at $p$. 
We say that $p$ is {\it unibranched}
if    $|\nu_C^{-1}(p)|=1$.
We denote by $C^{\nu}_p$ the partial normalization of $C$ at $p$;  the $\delta${\it-invariant}, $\delta_C(p)$,   of $p$  is $\delta_C(p)=p_a(C)-p_a(C^{\nu}_p).$  
 The following formula is well-known
\begin{equation}
 \label{eq-di}
  \delta_C(p)=\sum m_q(m_q-1)/2
\end{equation}
 where $q$ varies over all points infinitely near to $p$ (including $p$), and $m_q$ is the multiplicity of $q$. 
 
If $C,B\subset S$ are   reduced  curves   with no components in common, and $p\in C\cap B$, we write
$ 
(C\cdot B)_p 
$ 
for their multiplicity of intersection at $p$;
we say that $C$ and $B$ are   transverse  at $p$  if $(C\cdot B)_p=\mult_p(B)\mult_p(C)$.

We say that   $C$ is {\it hypertangent} to   $B$ if  $|\nu_C^{-1}(B\cap C)|=1$. 
 
 Let $S=\PP^2$, fix a reduced curve $B\subset \PP^2$  and a point
  $q\in B$.  The set of integral curves hypertangent to $B$ at $q$ is  
  $$
 \Hyp (B;q):=\{C\subset \PP^2:\ C \text{ integral}, \  C\cap B=\{q\},\  |\nu_C^{-1}(q)|=1\},
$$
and the set of all integral curves hypertangent to $B$ is
 $$\Hyp(B,1):=\cup_{q\in B} \Hyp(B;q).$$
For a  positive integer $d$ we     write
  $\Hyp_d(B;q)\subset  \Hyp (B;q)$ and  $\Hyp_d(B,1)\subset \Hyp(B,1)$, for the subsets parametrizing curves of degree $d$;
we view both of them as 
 subspaces of the projective space $\PP^{d(d+3)/2}$ parametrizing plane curves of degree $d$.

If $C\in  \Hyp_d(B;q)$ and $C$ is singular at $q$, then it necessarily has a unibranch singularity. We set  for $m\geq 1$
 $$
 \Hyp ^m(B;q)=\{C\in  \Hyp (B;q):  \  \mult_q(C)=m\}
$$
 and we define $\Hyp_d^m(B;q)\subset \Hyp_d(B;q)$ analogously. 

Now we extend to double intersections. As we said before,   $C$ is  {\it hyper-bitangent} to  $B$ if
$\nu_C^{-1}(C\cap B)\leq 2$ and 
we denote by
$ 
\Hyp(B,2) 
$ the set of curves hyper-bitangent to $B$.
We have
$\Hyp(B,1)\subset \Hyp(B,2)$.
We set
$$\Hyp_d(B,2):=\{C\in \Hyp(B,2):\  \deg C=d\}
$$
so that $\Hyp_d(B,2)\subset \PP^{d(d+3)/2}$.

 \subsection*{Acknowledgements.}  We are pleased to thank Laura Capuano, Wei Chen, Ciro Ciliberto, Pietro Corvaja, Kristin DeVleming, Edoardo Sernesi and  Umberto Zannier for discussions and useful comments about  this work. We are grateful to the referees for comments and suggestions which improved the paper.  LC is partially supported by  PRIN 2017SSNZAW and PRIN 2022L34E7W,  Moduli spaces and birational geometry.
 AT is partially supported by   PRIN  2022HPSNCR: Semiabelian varieties, Galois representations and related Diophantine problems and   PRIN   2020KKWT53: Curves, Ricci flat Varieties and their Interactions, and is a member of the INDAM group GNSAGA. 

\section{Hypertangency}
\label{sec-hyp} 
 
   \subsection{Preliminaries on unibranch points} We here state some basic facts about    unibranch points of curves; 
 we refer to     \cite{Wall} for an exhaustive treatise. 
   Let $q$ be a unibranch point of an integral curve $C\subset \PP^2$; 
   the tangent line  to $C$ at $q$   is  well defined,   we denote it by $L$.
When $L\neq C$ we say that
  $q$  is  an  $(m,n)${\it -point} if $\mult_q(C)=m$  and if  $ 
(C\cdot L)_q=n.
$ 

 If $m=1$ and $n\geq 3$  we say that  $q$ is a  {\it flex}  of $C$. 
 We   choose local   coordinates, $x,y$, so that $q=(0,0)$, the line $L$ has equation $y=0$, and 
$C$ has     equation
$f(x,y)=0$ with  $\deg f=d\geq n$  with
\begin{equation}
 \label{eq-f}
f(x,y)=a_{0, m}y^ m+   \sum_{ m+1\leq i+j\leq d} a_{i,j}x^iy^j 
\end{equation}
such that  $a_{0, m}\neq 0$ and   the smallest power of $x$   appearing in $f$  is  $x^n$, i.e. $a_{i,0}=0$ for $i<n$ and  $a_{n,0}\neq 0$.   
Notice that, as   will be clear in the sequel,   having local equation of type \eqref{eq-f} is not sufficient for $q$ to be   unibranched.

More generally, let   $S$ be any smooth  surface
and $q\in C\subset S$ a unibranch point of multiplicity $m$; let $D\subset S$ be a smooth  integral  curve,   $C\neq D$.  If $n:=(C\cdot D)_q>m$ we say that $D$ is tangent to $C$ at $q$, and  that 
 $q$ is an $(m,n)${\it -point of} $C$ {\it with respect to} $D$.

We are interested in the case where the surface $S$ is an   iterated  blow-up of $\PP^2$ over a fixed point $q^0\in\PP^2$, i.e.
we have a finite chain of blow-ups
$$S\la S^i\la \ldots \la S^0=\PP^2$$
all centered at a point   lying over $q^0$.  Let $C\subset S$ be an integral curve with a unibranch point $q\in C$. Let $L\subset S$ be the strict transform of a line  in  $\PP^2$ such that $q$ is an $(m,n)$-point with respect to $L$.
Consider the blow-up, $S'\to S$, at $q$,   let  $E\subset S'$ be the exceptional divisor and   $C'$  and $L'$ the strict transforms of $C$ and $L$. The map $\sigma:C'\to C$   induced by this blow-up is   bijective;    set $q':=\sigma^{-1}(q)$.
With this set-up, we have the following 
 \begin{lemma}
\label{lm-tan} 
 Let $0<m<n$ and  let $q\in C\subset S$  be  an  $(m,n)$-point with respect to $L$ 
 
 \begin{enumerate}[(a)]
 \item 
  \label{lm-tana}
If $n<2m$   then $q'$ is an
 $(n-m,m)$-point of $C'$ with respect to $E$.
   \item
    \label{lm-tanb}
  If $n>2m$  then $q'$ is an
 $(m,n-m)$-point of $C'$ with respect to $L'$.
  \item
     \label{lm-tanc}
     If $n=2m$ then $q'$ is an $m$-fold point of $C'$, and $C'$ is neither tangent to $E$ nor to $L'$ at $q'$.
 \end{enumerate}
\end{lemma}
\begin{proof} We can work locally and  use on $S$ affine coordinates so that the set-up described for \eqref{eq-f} holds.
To  blow-up at $q$ we set
  $y=vx$ and use $(x,v)$ as local coordinates in the blow-up at $q'=(0,0)$, the local equation of  the exceptional divisor $E$ is $x=0$, and   the local equation of $L'$ is  $v=0$.  Let $f'(x,v)=0$ be the affine equation of $C'$ obtained from \eqref{eq-f}:
$$
f'(x,v)=x^{-m}f(x,vx)=a_{0,m}v^m+\sum_{m+1\leq i+j\leq d} a_{i,j}x^{i+j-m}v^j.
$$
The smallest power of $x$ appearing  as a summand of $f'$ is $x^{n-m}$, and the smallest power of $v$ is $v^m$.
 We have  three cases.

  Case \eqref{lm-tana}  $n<2m$. Hence $n-m<m$ and we set $r=n-m$. Since $C'$ is unibranched at $q'$, and $f'$ contains the summand $x^r$ but not the summand $v^r$ (as $r< m$), 
all terms in $f'$ of degree at most $r$ divisible by $xv$ must vanish (for otherwise the lowest homogeneous part of $f'$ will be reducible, contradicting the fact that $q'$ is a unibranch point of $C'$).
  Hence the tangent line to $C'$ at $q'$ has local equation $x=0$,  hence  $C'$ is tangent  to $E$.  Since the smallest power of $v$   in $f'$ is $v^m$, we get that
 $q'$ is an $(r,m)$-point with respect to $E$.
 
  Case   \eqref{lm-tanb}  $n>2m$. 
The smallest power of $v$ appearing in $f'$ is $v^m$ and, arguing as in the previous case,
  $f'$ has no other term of degree at most $m$. Hence 
   the tangent line to $C'$ at $q'$ has local equation $v=0$,  so that  $C'$ is tangent  to $L$.  As the smallest power of $x$   is   $x^{n-m}$, we get that
 $q'$ is an $(m,n-m)$-point with respect to $L'$.   

  Case \eqref{lm-tanc}  $n=2m$. We   have both $x^{n-m}=x^m$  and   $v^m$ appearing in $f'$ as smallest powers of $x$ and $v$.
Since $C'$ is unibranched at $q'$, this is possible only if the homogeneous part of degree $m$ of $f'$,  has form
 $
f'_m=(\alpha v+\beta x)^m
$, 
with  $\alpha,\beta\neq 0$. Hence $q'$ is an $m$-fold point such that neither $E$ nor $L'$ are tangent to $C$ at $q'$.  \end{proof}
  We will frequently use the following well known facts.
\begin{remark}     
\label{rk-int}
Let $B,C\subset S$ be two  integral   curves  and $q\in B\cap C$.
Suppose that  $B$ is smooth at $q$, and that   $C$  has a unibranch 
 $m_C$-fold  point at $q$. 
  If $S'\to S$ is the blow up at $q$ and $C', B'\subset S'$ are the strict transforms of $C$ and $B$,
  then
  $ 
(C'\cdot B') =(C\cdot B) - m_C.
$ 

Moreover,
 $(C\cdot B)_q>m_C$ if and only if $C$ and $B$ are   tangent   at $q$.
     \end{remark}
 
    \subsection{Hypertangency at unibranch points of plane curves}
   We   describe plane   curves  hypertangent to  a  curve  at a smooth point. 
  See Subsection~\ref{ss-ex} for some examples.
     Recall that a point $q$ of a plane curve $C\subset \PP^2$ of degree at least $2$ is   an $(m,n)$-point if it is unibranch, $m=\mult_q(C)$ and
    $(L\cdot C)_q=n$, where $L$ is the tangent line to $C$ at $q$.

\begin{thm}
\label{mirror-sm}
 Let $B, C\subset \PP^2$ be   integral curves  of degree at least $2$ such that 
 $B\cap C=\{q\}$; assume that   $q$ is a $(1,l)$-point    for $B$ and  an $(m,n)$-point for $C$.  
Then   $n=lm$ and 
  the $\delta$-invariant of $q$ on $C$ satisfies
$$
\delta_C(q)\geq (m-1)(\deg B\deg C-m)/2.
$$

\end{thm}

\begin{proof}
 Set $b=\deg B$, $d=\deg C$,    and let $q$ be an $(m,n)$ point of $C$.
 Hence $l\leq b$ and $1\leq m<n\leq d$.
  By hypothesis,  $$(C\cdot B)_q=(C\cdot B)=bd >m$$
 therefore $C$ and $B$    must have the same tangent line at $q$ 
(by Remark~\ref{rk-int}); we write $L$ for this line.
Let $C^1$,  $B^1$  and $L^1$ be the strict transforms of $C, B$ and $L$ in the blow-up, $S^1$, of $\PP^2$ at $q$,
and let $q^1\in C^1$ be the   point lying over $q$. By hypothesis $C^1$ and $B^1$   meet only in $q^1$, hence    (as $d>m$, $b\geq 2$)
$$
(C^1\cdot B^1)_{q^1}=(C^1\cdot B^1)=bd-m>bm-m\geq 2m-m=m.
$$
Now, $B^1$ is smooth at $q^1$ and  $\mult_{q^1}(C^1)\leq \mult_q(C)=m$, hence
 $$(C^1\cdot B^1)_{q^1}>m\geq \mult_{q^1}(C^1)\mult_{q^1}(B^1),$$ hence $C^1$ and $B^1$   are  tangent   in $q^1$. 

We   prove that $C$ has an $(m,lm)$-point at $q$.
Suppose $l=2$;
   Lemma~\ref{lm-tan} implies that $B^1$ is neither tangent to $L^1$, nor to the exceptional divisor, hence the same holds for $C^1$. Hence, by the same Lemma, $q$ is a $(m,2m)$ point for $C$, and we are done.
 
  To set-up an inductive argument  
we introduce a sequence of blow-ups as follows. We already considered the blow-up, $S^1$, of $S^0=\PP^2$ at $q^0=q$. For $i\geq 1$ we let $S^i$ be the blow-up of $S^{i-1}$ at the unique point $q^{i-1}\in C^{i-1}$
lying over $q$. We denote by $C^i$, $B^i$  and  $L^i$ the strict transforms of $C,B$ and $L$.
Notice   that  the multiplicity of $C^i$ at $q^i$ is   at most $m$.
 
\begin{claim}
Let $l\geq 3$; 
for every $1 \leq i \leq  l-2$ the point $q^i$ is an $(m, n-im)$-point with respect to $L^i$ for  $C^i$, 
moreover
$C^i$ and $B^i$ are tangent   at 
$q^i$.   
\end{claim}

We prove the claim by induction.
If $i=1$ then  $q$ is a $(1,l\geq 3)$-point for $B$, hence Lemma~\ref{lm-tan} gives that $B^1$ is   tangent to $L^1$ over $q$. We proved earlier that $B^1$ and $C^1$ are tangent in $q^1$, hence    $C^1$ is tangent to $L^1$ at $q^1$; hence, by Lemma~\ref{lm-tan}, $q$ is an $(m,n)$-point for $C$ with $n>2m$, and $q^1\in C^1$ is  an
$(m, n-m)$-point with respect to $L^1$. The proof of the  base   is complete.

Suppose $C^{i-1}$ has a $(m,m-(i-1)n)$-point with respect to $L^{i-1}$ at $q^{i-1}$,
 and $B^{i-1}$ and  $C^{i-1}$ are tangent to $L^{i-1}$ at $q^{i-1}$. Now $ B^{i-1}$ has a
$(1,l-(i-1))$-point  in  $q^{i-1}$; as $i\leq l-2$ we have  $l-(i-1)\geq l-l+3=3$, hence $B^i$ is tangent to $L^i$.
Now, $l\leq b$, hence    $i\leq l-2\leq b-2$. Therefore 
$$
(C^i\cdot B^i)_{q^i}=(C^i\cdot B^i)=bd-im>bm-im\geq bm-(b-2)m=2m>m. 
$$
Hence $C^i$ and $B^i$ are tangent in $q^i$ at $L^i$. Therefore $C^i$ is tangent to $L^i$ and  case \eqref{lm-tanb} of Lemma~\ref{lm-tan}    must occur for $C^{i-1}$, i.e. $m-(i-1)n>2m$  and $C^i$ has a
 $(m,n-im)$-point with respect to $L^i$.
 The     claim is proved.
 
 \
 
 Thus $C^{l-2}$ has an $(m,n-(l-2)m)$-point with respect to $L^{l-2}$ at $q^{l-2}$ , and   $B^{l-2}$ and $ C^{l-2}$
  are   tangent to $L^{l-2}$ at  $q^{l-2}$. Now,  $B^{l-2}$ has a $(1,2)$-point at 
 $q^{l-2}$, hence its strict transform,  $B^{l-1}$, is neither tangent to  $L^{l-1}$ nor to the exceptional divisor.
 We have,  as $m<d$ and $l\leq b$  
 $$
 (C^{l-1}\cdot B^{l-1})_{q^{l-1}}=(C^{l-1}\cdot B^{l-1})=bd-(l-1)m>bm-(b-1)m=m,  
 $$
hence  $C^{l-1}$ and $ B^{l-1}$ are tangent in $q^{l-1}$, hence $C^{l-1}$ is neither tangent to  $L^{l-1}$ nor to the exceptional divisor. Now case \eqref{lm-tanc} of Lemma~\ref{lm-tan}  occurs for $C^{l-2}$, i.e. 
  $q^{l-2}$ is a $(m,2m)$-point. Therefore $2m=n-(l-2)m$, hence $n=lm$, as stated.

Let us now study the $\delta$-invariant for $q$ as a point of $C$.
Set
 $$
 h=\lceil (bd-m)/m\rceil.
 $$
We   prove, by induction on $i$, that $q^i$ is an $m$-fold point of $C^i$   for every $i\leq h-1$,
and   $B^i$ and $C^i$ are tangent in $q^i$. 
We already proved this  for every $i\leq l-1$, hence the base case is settled and we   assume $i\geq l$. 
The argument is similar to the one used in the previous part.
  Assume $\mult_{q^{i-1}}(C^{i-1})=m$, and $B^{i-1}$  tangent to  $C^{i-1}$    in $q^{i-1}$. 
Then, as $i\leq h-1$, we have
 \begin{align*}
  (C^{i}\cdot B^{i})_{q^{i}}= (C^{i}\cdot B^{i})=bd-im\geq bd-(h-1)m=  \\   
  bd-(\lceil (bd-m)/m\rceil-1)m 
> bd-( bd/m -1)m=bd-bd+ m=m
\end{align*}
(as $\lceil (bd-m)/m\rceil<bd/m$).
Hence $(C^{i}\cdot B^{i})_{q^{i}}>m$, hence $C^i$ and $B^i$ are tangent at $q^i$. Now,   the curve $B^{i-1}$ has a  smooth point at $q^{i-1}$, hence $B^i$ is not tangent to the exceptional divisor in $q^i$.
Hence the same holds for $C^i$. Hence case \eqref{lm-tana} of Lemma~\ref{lm-tan}  does not occur for $C^{i-1}$,  hence $q^i$ is an $m$-fold point of $C^i$. So we are done.
  
  Since $C^i$ has an $m$-fold point at $q^i$ for every $i=0,\ldots, h-1$, by \eqref{eq-di} the $\delta$-invariant of $q$ satisfies
 \begin{align*}
  \delta_C(q)&\geq hm(m-1)/2\\
 &= \lceil (bd-m)/m\rceil m(m-1)/2  \\
& \geq ((bd-m)/m )m(m-1)/2\\&=   (bd-m) (m-1)/2.
\end{align*}
 The proof is now complete.
 \end{proof}

   \section{Hyper-bitangency for     3C-curves}
   \label{sec-3C} 

  \subsection{Definition and   simple cases.} We  
 study     hyper-bitangent curves to a curve $B\subset \PP^2$ which   is the transverse union of three integral curves.  
\begin{defi}
\label{defi-3C}
 A {\it 3C-curve} is a reduced plane curve   $B=  B_1\cup B_2\cup B_3$,   
with $B_i$  integral   of degree $b_i\geq 1$,   such that  every point in $B_i\cap B_j$ is a node of $B$ for all $i\neq j$.
We always assume $b_1\leq b_2\leq b_3$.
We set
 $$B_i\cap B_j=\{p_{i,j}^{t},\quad t=1,\ldots, b_ib_j\} 
 $$
with $p^t_{i,j}=p^t_{j,i}$; we often omit the superscript $t$.
Notice that the components of $B$ meet only pairwise, and  transversally. 
 We write  $N:=\cup_{i\neq j} B_i\cap B_j.$   
 \end{defi}

  We begin with the case $b_1=b_2=b_3=1$.
  This is a particularly simple curve  which can be     easily handled.
  \begin{prop}\label{triangle}
 Let  $B= B_1\cup B_2\cup B_3$   be a 3C-curve of degree  $3$.
 Then $\dim \Hyp_d(B,2)\geq 1$ for every $d\geq 1$.
 \end{prop}
\begin{proof}
The curve $B$ is the union of three lines;  set $B_i\cap B_j=\{p_{i,j}\}$ so that $N=\{p_{1,2},p_{1,3},p_{2,3},\}$. 
 
 Suppose $d=1$; then the one-dimensional space of  lines through $\{p_{i,j}\}$,  with  $B_i$ and $B_j$ removed,  lies in $\Hyp_1(B,2)$, and these are the only elements of $\Hyp_1(B,2)$. Hence  $\dim \Hyp_1(B,2)=1$.
 
Let $d\geq 2$ and $C\in \Hyp_d(B,2)$. Obviously $C\cap N\neq \emptyset$, say $p_{1,2}\in C$. Then the tangent line to $C$ at $p_{1,2}$ (which is well defined as $C$ is unibranched at $p_{1,2}$) must be different from at least one of the (different) tangent lines to  $B_1$ and $B_2$ at $p_{1,2}$. Hence $C$
must be transverse to 
 at  least one between $B_1$ and $B_2$, say $C$ transverse to $B_1$, hence $C$ meets $B_1$ in a further point. Since $C$ must also meet $B_3$ we get  $p_{1,3}\in C$, and $C$ cannot be transverse to $B_3$.  We derive that $C$ is hypertangent to $B_2$ and $B_3$ respectively at  $p_{1,2}$ and  $p_{1,3}$. Now,
setting 
 $m =m_{p_{1,2} }(C)$ and  $n =m_{p_{1,3} }(C)$, we have
  $$
 d=  (C\cdot B_1)=  (C\cdot B_1)_{p_{1,2}}+ (C\cdot B_1)_{p_{1,3}}=m+n\leq d
$$
hence   $d=n+m$.   Interchanging the three components, we   derive  
$$   \Hyp_d(B,2) =  
  \bigcup_ {m=1}^{d-1}  \bigcup _{\stackrel {\quad i,j,h=1,2,3}{ i\neq j}}  \Hyp_d^m(B_i;p_{h,i})\cap \Hyp_d^{d-m}(B_j;p_{h,j}).
 $$
Consider  $m=1$ and the subspace 
$ \Hyp_d^1(B_i;p_{h,i})\cap \Hyp_d^{d-1}(B_j;p_{h,j}).$
Now, $ \Hyp_d^1(B_i;p_{h,i})$ is the space of integral degree-$d$ curves    passing through     $p_{h,i}$, 
and meeting the line $B_i$ with multiplicity $d$ at $p_{h,i}$. It is easy to see that the closure of this space in 
$\PP^{d(d+3)/2}$ is a linear subspace of codimension $d$.
 
Next,  $ \Hyp_d^{d-1}(B_j;p_{h,j})$  is the space of integral degree-$d$ curves  having a $(d-1,d)$-point at 
$p_{h,j}$ with tangent line equal to $B_j$. To compute the codimension of its closure in $\PP^{d(d+3)/2}$
 we can assume that $p_{h,j}$ is the origin and the line $B_j$ has equation $y=0$. Then
a curve $C\in  \overline{\Hyp_d^{d-1}(B_j;p_{h,j})}$ has   equation $\sum_{d-1\leq i+j\leq d}a_{i,j}x^iy^j=0$ with $a_{i,j}=0$
for every $i+j=d-1$ and $i\neq 0$. One easily checks that such polynomials form a linear subspace of codimension equal to  $d(d-1)/2+d-1=(d^2+d-2)/2$. 
Therefore
$$
\dim  \overline{\Hyp_d^1(B_i;p_{h,i})}\cap \overline{\Hyp_d^{d-1}(B_j;p_{h,j})}\geq  d(d+3)/2-d-(d^2+d-2)/2 =1.
$$
To conclude the proof it suffices to prove that  a general point of the above intersection parametrizes an integral curve, i.e. an element in 
$\Hyp_d(B,2)$.   To do this, it suffices to prove that the above intersection contains one integral curve,
i.e. that there exists an integral curve of degree $d$ meeting  $B_i$ with multiplicity $d$ at $p_{h,i}$ and with a $(d-1,d)$-point at 
$p_{h,j}$ with $B_j$ as tangent line.   By choosing projective coordinates $x,y,z$ such that 
 $B_i$ and $B_j$  have respective equations $z=0$ and   $y=0$, with $p_{h,i}=(0:0:1)$ and $p_{h,j}=(0:1:0)$, we have the curve  of equation $y^{d-1}z=x^d$ (see subsection~\ref{sec-snc} for the properties of such a curve).
\end{proof}

  \subsection{Hyper-bitangent lines}
We now describe   hyper-bitangent lines to a 3C-curve of   degree at least $4$.
This is  quite elementary, possibly part  of it already  known. We include it for completeness and lack of references.
\begin{prop}\label{prop-3C1}
 Let $B= B_1\cup B_2\cup B_3$ be a  3C-curve of degree $b\geq 4$. 
\begin{enumerate}[(a)]
  \item
  \label{prop-3C1a}
If  $b=4$ then $|\Hyp_1(B,2)|=10$. 

More precisely, $C\in  \Hyp_1(B,2)$ if and only if $C$ is one of the four lines through $B_3\cap(B_1\cup B_2)$ different from $B_1$ and $B_2$, or $C$ is one of the six lines tangent to $B_3$ and passing through a point in $N = \cup_{i \neq j} B_i\cap B_j$.  
   
\begin{center}
\begin{tikzpicture}[xscale=1.1,yscale=1.1]
 \draw[](0.8,1) ellipse (1.5 and 0.9) node at (0,0){$B_3$}; 
  \draw [-] (0.5,2.8)--(2,-0.2)node at (2.2,0){$B_2$};
    \draw  [-] (-1,-0.2)--(1,2.8)node at (-1.1,0){$B_1$};
\draw[thick, red] (1,2.6)--(-1.8,0.5);
   \node at (-1.2,1.2){\color{red}$ C$};
   \draw[thick, blue] (-1, 0.58)--(2.3,0.27);
     \node at (-1.2, 0.55){\color{blue}$ C$};
 \draw[thick,purple] (0.3, 2 )--(3,1.6);
  \node at (2.9,1.8){\color{purple}$ C$};
\end{tikzpicture}

  \end{center}
   \item
     \label{prop-3C1b}
If $b\geq 5$  then $\Hyp_1(B,2)$ is finite, and it  is empty  if   $B$ is general.

More precisely, if $b_1=b_2=1$ then $|\Hyp_1(B,2)|\leq 3b_3(b_3-2)$, and $|\Hyp_1(B,2)|\leq 2|N|$ otherwise  (recall  $N=\cup_{i\neq j} B_i\cap B_j$).

  \end{enumerate}
\end{prop} 
 \begin{proof}
 If $b=4$ then   $b_1=b_2=1$ and $b_3=2$;  set  $B_1\cap B_2=\{p_{1,2}\}$. We first look at the lines through two points of $N$. If $C$ is a line  through $p_{1,3}^{t}$ and $p_{2,3}^{t'}$ 
it clearly lies in $\Hyp_1(B,2)$. If $C$ is a line  through $p_{i,3}^{t}$ and $p_{i,3}^{t'}$  then it is equal to $B_i$ which is not possible.  Now suppose   that $C\in \Hyp_1(B,2)$ is not one of these   lines; then we must have $p_{1,2}\in C$  and $C$ must meet $B_3$ in a unique point, hence $C$ must be  tangent to $B_3$. Part \eqref{prop-3C1a} is proved.

Let  $b\geq 5$ and let   $B$  a general curve; we  can make  the following assumptions. 
If $b_3\geq 3$  then  $B_3$ has finitely many  flexes, hence   finitely many  flex lines, and finitely many bitangent lines;  we assume that $B_1$ and $B_2$ do not pass through any flex  of $B_3$
 or through any point where $B_3$ meets its bitangent lines, and  that $B_1\cap B_2$ intersects no  line   hyper-bitangent to $B_3$. 

 If $b_3=2$ then $b_2=2$,
 we  assume that $B_1$ does not intersect $B_2\cup B_3$ at any point where  $B_2\cup B_3$  meets its  bitangent lines. 
 
 Finally,  given  any  $p_{2,3}\in B_2\cap B_3$, there are   finitely many lines through $p_{2,3}$ that are tangent to $B_3$, hence there are   finitely many points  $r\in B_3$ such that the tangent line to $B_3$ at $r$ passes through $B_2\cap B_3$; we   assume that $B_1$ does not intersect $B_3$ in any of these points. Also, if $b_2\neq 1$,  we assume that $B_1$ does not intersect $B_2$ in any   point lying on some tangent line to $B_2$ or $B_3$ passing through     $p_{2,3}$.

By contradiction, let $C\in  \Hyp_1(B,2)$.
Suppose  $|C\cap N|=2$. 

If  $b_3\geq 3$  then, as $N$ contains no flex of $B_3$, we have  $C\cap B=\{p_{1,3},p_{2,3}\}$ and,
as $C$ cannot   be   a bitangent to $B_3$,  we have
$$
b_3=(C\cdot B_3)_{p_{1,3}}+(C\cdot B_3)_{p_{2,3}}\leq 3
$$ 
hence $b_3=3$.
Now $C$ must be tangent to $B_3$ in one of the points, $p_{i,3}$. 
Hence $C$ is a line through $B_2\cap B_3$ tangent to $B_3$ and intersecting $B_1\cap B_3$; we excluded the existence of such lines, so we are done.

Let $b_3=2$, hence $b_2=2$.
 If $b_1=1$,  up to switching $B_2$ and $B_3$ we have  only the case $C\cap B=\{p_{1,2},p_{2,3}\}$.
Then  $C$ is tangent to $B_3$ at $p_{2,3}$ and intersects $B_1\cap B_2$,   which is excluded.

Let  $b_1= 2$.  Our generality assumptions prevent us from switching $B_1$ with $B_2$ or $B_3$, so we have more cases. If $C\cap B=\{p_{1,2},p_{1,3}\}$  then $C$ is  a bitangent of $B_2\cup B_3$, which is excluded (as before).   If $C\cap B=\{p_{1,2},p_{2,3}\}$   then  $C$ is tangent to $B_3$ at $p_{2,3}$ and intersects $B_1\cap B_2$, which is   excluded. As we can switch $B_2$ with $B_3$ we are done.
We thus proved that $|C\cap N|=1$.
 
 Suppose
  $C\cap N=\{p_{1,2}\}$.  
Then $C$  meets $B_3$ in a point $r\not\in N$, and it is hypertangent to $B_3$ at $r$.
If $b_3\geq 3$ then $r$ is a flex and, by our assumptions, the flex line does not pass through $B_1\cap B_2$.  If $b_3=2$ then $b_2=2$ and $C$ is a bitangent of $B_2\cup B_3$, which is also excluded.

Suppose $C\cap N=\{p_{1,3}\}$, then either $b_3\geq 3$ and $C$ is a flex line of $B_3$ which is excluded, or  $b_3=b_3=2$ and    $C$ is a bitangent of $B_2\cup B_3$, which is excluded.

Suppose $C\cap N=\{p_{2,3}\}$. 
Now, $C$ cannot be tangent to both $B_2$ and $B_3$,
hence $b_2=1$, hence $b_1=1$, hence $b_3\geq 3$, hence   $p_{2,3}$ is a flex of $B_3$. A contradiction.
We thus proved that $\Hyp_1(B,2)$ is empty for  $B$ general.

If  $B$ is an arbitrary curve, the proof   shows that for  a curve $C\in \Hyp_1(B,2)$ only  two cases can occur.
First  case: $C$ is
tangent to some component of $B$ in a point of $N$; since at each point of $N$ there are   two  such tangent lines  we have at most $2|N|$  possibilities for such a $C$.

Second case:  $C$ meets all components of $B$ transversally along $N$. Then one easily checks that $b_1=b_2=1$, moreover $C$ passes through $B_1\cap B_2$ and is hypertangent to  $B_3$ in a flex (or rather a hyper-flex) or in a unibranch singular point. 
It is well known that the number of flexes of $B_3$ is at most equal to  $3b_3(b_3-2)$;
hence in the second case we have at  most $3b_3(b_3-2)$ possibilities for $C$.

If
 $b_1=b_2=1$  the first case occurs only 
with $C$ hypertangent  $B_3$, hence the bound $ \Hyp_1(B,2)\leq 3b_3(b_3-2)$ holds.
 \end{proof}
   \subsection{Hyper-bitangent curves of higher degree}
 We begin with a geometric description of hyper-bitangent curves.
\begin{thm}\label{3C}
Let $d\geq 2$.
 Let  $B= B_1\cup B_2\cup B_3$   be a 3C-curve of degree  $b\geq 4$ such that
 $\Hyp_d(B,2)$ is not empty. 
   Then $b_1=1$  and   the following  occur.
   \begin{enumerate}[(a)]
 \item
 \label{3Ca}
If $b_2= 1$, setting $B_1\cap B_2=\{p\}$, we have
$$\Hyp_d(B,2) =
 \bigcup_{\stackrel {i=1,2}{q\in B_i\cap B_3}}  \Hyp_d^{d-1}(B_i;p)\cap \Hyp_d (B_3;q).
 $$
 \begin{center}
\begin{tikzpicture}[xscale=1.2,yscale=1.2]
\draw[blue](-0.5,2)to  [out=  0, in= 100]  (1.5,1.36);
\draw[blue](1.5,1.36)to  [out=  270, in= 20]  (-0.5,1.2);
\draw[blue](-0.5,1.2)to  [out=  190, in= 150]  (-0.7,0.8);
\draw[blue](-0.7,0.8)to  [out= 320, in= 150]  (2,0.2)node at (0.8,0.25){$B_3$};

  \draw [blue] (0.5,2.8)--(2,-0.2)node at (2.2,0){$B_2$};
 \draw  [blue] (-0.9,-0.2)--(1,2.8)node at (-1,0){$B_1$};
 
    \node at (0.90,2.35){\tiny$p$};
        \node at (-0.13,0.6){\tiny$q$};
\draw[thick, red](0.6,3)to  [out=  -100, in= 100]  (0.72,2.36);
\draw[thick, red](0.2,2.8)to  [out= -10, in= 120]  (0.72,2.36);
 \draw[thick, red](0.2,2.8)to  [out= 160, in= 170]  (-1,0.8);
\draw[thick, red](-1,0.8)to  [out=   -10, in= 160]  (-0.4,0.6);
  \draw[thick, red](-0.4,0.6)  to  [out= -30, in= 170]  (0.5,0.1);
    \node at (-0.3,2.5){\color{red}$ C$};
    \end{tikzpicture}
  \end{center}

  \item
 \label{3Cb}
If $b_2\geq 2$ then $\Hyp_d(B,2)$ is empty for $d\geq 3$ and  $$
\Hyp_2(B,2) = \bigcup_{\stackrel  {p\in B_1\cap B_2}{q\in B_1\cap B_3} } \Hyp_2 (B_2,p)\cap \Hyp_2 (B_3,q).$$
  
    \item
 \label{3Cc}
Every $C\in \Hyp_d(B,2)$ is rational (hence  $\cE (B)=\Hyp (B,2)$, with $\cE(B)$    defined in \eqref{eq-EB}).
   \end{enumerate}
    
\end{thm}
 \begin{proof}
Since $b\geq 4$ we have $b_3\geq 2$.  Let $C\in \Hyp_d(B,2)$,
then  $|C\cap B|\leq 2$; on the other hand $B_1\cap B_2\cap B_3=\emptyset$, hence
$|C\cap B|=2$. 

Let us prove that $C\cap B_3\subset N$ and  $|C\cap B_3|=1$.
  By contradiction, suppose $C\cap B_3$ contains a point not in $N$. Hence $C$ must intersect $B_1\cup B_2$ in exactly one point,    $p_{1,2}\in B_1\cap B_2$. Now,  $p_{1,2}$ is  a unibranch $n$-fold point  of $C$, for some $n <d$  (as $d>1$). 
Since $B_1$ and $B_2$ meet  transversally, $C$ is tangent to one of them and transverse to the other,  say  $C$ is transverse to $B_i$ with $i<3$. Therefore  $$b_id= (C\cdot  B_i)=  (C\cdot B_i)_{p_{1,2}} =n<d$$
 a contradiction. Hence $C\cap B_3\subset N$.
 
By contradiction, suppose $|C\cap B_3|=2$;
as $C$ must intersect   $B_1$ and $B_2$ we have $C\cap B_3=\{p_{1,3},p_{2,3}\}\subset N$, with $p_{i,3}\in B_i\cap B_3$. 
Since $B_3$ and $B_i$ meet transversally, $C$ must be transverse to either $B_i$ or $B_3$ at $p_{i,3}$. If $C$ is transverse to $B_i$ then, arguing as above, we get $(C\cdot  B_i)<d$, a contradiction.
Hence  $C$ is transverse to $B_3$ at both $p_{1,3}$ and $p_{2,3}$.  Set $m_i=\mult_{p_{i,3}}(C)$, so that  $m_1+m_2\leq d$. Then 
$$
b_3d= (C\cdot B_3)=  (C\cdot B_3)_{p_{1,3}}+  (C\cdot B_3)_{p_{2,3}}=m_1+m_2 \leq d
$$
which is impossible, as $b_3\geq 2$. 

We thus proved that   $C\cap B_3=\{p_{i,3}\}$   for one   $i\in \{1,2\}$. 
Hence
$C$ is hypertangent to $B_3$ at $p_{i,3}$.
As $C$ is not transverse to $B_3$ a $p_{i,3}$,
 it must be transverse to $B_i$, hence it must meet $B_i$ in a further point. As $C$ must   meet the other component, $B_j$ with $j\neq i,3$, there exists a point $p_{1,2}\in B_1\cap B_2$ such that  $p_{1,2}\in C$. We obtain
$C\cap B=\{p_{1,2},p_{i,3}\}$  and  $C\cap B_j=\{p_{1,2}\}$. 
Now $B_j$ and $C$ meet only at $p_{1,2}$, hence
 $C$  cannot be transverse to $B_j$ at  this point, hence it must be transverse to $B_i$. Therefore
 $$
b_id=  (C\cdot B_i)=  (C\cdot B_i)_{p_{1,2}}+ (C\cdot B_i)_{p_{i,3}}=\mult _{p_{1,2}}(C)+\mult_{p_{i,3}}(C)\leq d
$$
hence $b_i=1$ and equality   holds, i.e.
$  \mult _{p_{1,2}}(C)+\mult_{p_{i,3}}(C)=d.$ 
 Therefore $b_1=1$ and
 $C\in \Hyp_d^{d-m}(B_2,p_{1,2} )\cap \Hyp_d ^{m}(B_3, p_{1,3} )$, 
 with $m=\mult _{p_{1,3}}(C)$.
 
We now show that $m=1$. Since $B_3$ has degree at least $2$ and is smooth at $p_{1,3}$,  it has a $(1,l)$-point there for some $l\geq 2$. We can apply Theorem~\ref{mirror-sm} to $B_3$ and $C$,  getting that $p_{1,3}$ is an $(m,lm)$ point for $C$. Suppose $m\geq 2$;
  the same theorem yields  $\delta_C(p_{1,3})\geq (b_3d-m)(m-1)/2$. On the other hand 
  $C$ has also a $(d-m)$-fold point at  $p_{1,2}$, hence   $\delta_C(p_{1,2}) \geq (d-m)(d-m-1)/2$.
  Therefore  
  
\begin{align*}
 g(C)&\leq {\binom{d-1}{2}}-(b_3d-m)(m-1)/2-(d-m)(d-m-1)/2\\
& ={\binom{d-1}{2}}-(d^2-d-2md+b_3md-b_3d+2m)/2\\
 &= {\binom{d-1}{2}}-(d^2 +d(-1+  m(b_3-2)-b_3)+2m)/2\\
  &\leq   {\binom{d-1}{2}}-(d^2-3d+4)/2\\
   &= (d^2-3d+2)/2-(d^2-3d+4)/2\\
   & <0,
\end{align*}
 as $m\geq 2$ and $ b_3\geq 2$. This is impossible. Hence $m=1$ and 
 \begin{equation}
\label{m23}
  C\in \Hyp_d^{d-1}(B_2,p_{1,2} )\cap \Hyp_d   (B_3, p_{1,3} ). \end{equation} 
  If $b_2=1$ we can switch roles between $B_1$ and $B_2$; part \eqref{3Ca} is proved.
 
Assume  $b_2\geq 2$.  Then, as  $B_2$   is smooth at $p_{1,2}$,  we can apply Theorem~\ref{mirror-sm},
which gives that $C$ has an $(h,lh)$ point at $p_{1,2}$, for some $h\geq 1$. Now  \eqref{m23} implies that
$p_{1,2}$  is  a $(d-1,d)$-point of $C$. Therefore
$(h,lh)=(d-1,d)$, hence $h=1$ and $d=2$. Therefore $\Hyp_d(B,2)$ is empty if $d\geq 3$, and the proof of part \eqref{3Cb} is complete.

A curve of degree $d$ having a $(d-1)$-fold  point  is necessarily rational, hence part \eqref{3Cc} follows  
 from the previous parts.   \end{proof}
 
 If $d\geq 2$, Theorem~\ref{3C} implies that  $\Hyp_d(B,2)=\emptyset $ whenever  $b_1>1$, 
  or $b_2\geq 2$ and $d\geq 3$.
  We now treat  the remaining cases.  
We denote by $\Hyp_{\geq 2}(B,2)$ the set of   curves of degree at least $2$ that are hyper-bitangent to $B$, i.e. 
$\Hyp_{\geq 2}(B,2):=\cup_{d\geq 2}\Hyp_{d}(B,2)$.

\begin{thm}\label{main}
Let $B$ be a  3C-curve   of degree  $b\geq 4$
such that  $b_1=1$.  Then  \begin{enumerate}[(a)]
        \item  \label{maina}
  $|\Hyp_{\geq 2}(B,2)|\leq b_3 \max\{2,b_2\} $;  
       \item
 \label{mainc}
 $|\Hyp_{\geq 2}(B,2)|=0$  if  $B$ is general.
\end{enumerate}
\end{thm}

\begin{proof} 
 Let $C\in \Hyp_d(B,2)$ with $d\geq 2$;
by Theorem~\ref{3C}, up to switching $B_1$ and $B_2$ when $b_2=1$, we have  $C\in   \Hyp_d^{d-1}(B_2;p)  \cap \Hyp_d  (B_3;q)$ with $p\in B_1\cap B_2$ and $q\in B_1\cap B_3. $
Now,    $C$ has a $(d-1,d)$-fold point at $p$ where it is tangent to  $B_2$, and a smooth point at $q$ where it is tangent to $B_3$. We can choose homogeneous coordinates $(X,Y,Z)$ in $\PP^2$ so that  $p=(0:1:0)$ and the tangent line to $C$ at $p$ has equation $z=0$. 
Therefore in the open subset where $Y\neq 0$ the curve $C$ has affine  equation 
$ 
z^{d-1} 
=\sum_{i=0}^dc_ix^iz^{d-i}.
$ 
 We  can assume that $q=(0:0:1)$ and the tangent line to $B_3$ has equation $y=0$. Hence $c_0=c_1=0$ and the affine equation of $C$  where  $Z\neq 0$ is
$$
y=g(x) \quad \quad\text{ where }\quad \quad g(x):=c_dx^d+c_{d-1}x^{d-1} +\ldots +c_2x^2.
$$
\begin{claim}
  $g(x)=c_dx^d$.
\end{claim}
We can  assume $d\geq 3$. We provide two  proofs of the claim, giving  different insights.

First proof. We follow the proof of \cite[Lm. 42]{Kollar}. The tangent line, $L$, to $C$ at $q$ has equation $y=0$, hence it suffices to prove that $(L\cdot C)_q=d$. 
Denote by $(B_3)_C$ and $L_C$ the (Cartier) divisors cut by $B_3$ and  $L$ on $C$. Set  $\n:=b_3$ to simplify; we have $(B_3)_C=ndq$ and, of course, $(B_3)_C \sim  nL_C$, hence $n(dq-L_C)\sim 0$. On the other hand, $C$ is a rational curve whose unique  singular point, $p$, is unibranch, hence $\Pic C$ has no torsion; see \cite[Ex. 6.11.4]{Hartshorne}. Therefore $dq\sim  L_C$,  hence (as  curves of degree  $d$ cut on $C$ a complete linear series) $dq= L_C$, and the claim is proved.

Second proof.
Let $f(x,y)=0$ be the affine equation of $B_3$;  as $B_3$ is smooth at $q$ and tangent to  $y=0$, we have
\begin{equation}
 \label{eq-B3}
f(x,y)=y+\sum_{2\leq i+j\leq \n }a_{i,j}x^iy^j
\end{equation}
with $\n=b_3$.
We have $(C\cdot B_3)=(C\cdot B_3)_q=d\n $, therefore  $
f(x,g(x))=\lambda x^{d\n}
$ 
for some $\lambda\neq 0$. Hence the following   holds
\begin{equation}
 \label{eq-dn}
g(x)+\sum_{2\leq i+j\leq \n }a_{i,j}x^ig(x)^j=\lambda x^{d\n}.
\end{equation}
Now, for every $i,j$ as above, the product $x^ig(x)^j$ is a sum of monomials in $x$ whose degrees range in  the set $I(i,j) $, where
$$
I(i,j) =[i+2j,i+jd]\cap \N.
$$
We have
$$
i+dj \leq  d(i+j)\leq d\n 
$$
with equality  if and only if  $j=\n$ and $i=0$.
 We have 
 $$
 I(0,\n) =[2\n, d\n ], \quad  \quad  I(1,\n -1)=[2\n -1,d\n  -d+1],
 $$
therefore  in the left   side of \eqref{eq-dn}, the following monomials  (up to scalar)
 $$
 x^{d\n }, \  x^{d\n -1}, \ldots, x^{d\n -d+2}
 $$
 appear  only in $g(x)^{\n}$. For   \eqref{eq-dn} to hold, in its  left  side  

(a) the coefficient of $ x^{d\n }$ is non zero, hence $a_{0,\n}\neq 0$ and $c_d\neq 0$;

(b)  the coefficients of $ x^{d\n -1}, \ldots, x^{ d\n-(d-2)}$  are zero.

\noindent
We now show, by   induction on $h$, that (a) and (b) imply   that $c_{d-h}=0$ for $h=1,\ldots, d-2$.  
 We write $g(x)^\n$ as follows
$$
g(x)^n=\sum_{k=1}^{d-2} \Bigg(\sum_{\stackrel{2\leq i_1<\ldots<i_k\leq d}{\sum_1^k n_j=\n ,\ n_j\geq 1}}
\mu_{n_1,\ldots, n_k}c_{i_1}^{n_1}\cdot\ldots\cdot c_{i_k}^{n_k} x^{n_1i_1+\ldots + n_ki_k}\Bigg)
$$
 where $\mu_{n_1,\ldots, n_k}$ are positive integers which we can ignore.
Since $i_k\leq d$ and $i_j<i_{j+1}$ we have   $i_j\leq i_k-(k-j)\leq d-(k-j)$, hence the exponent  of $x$ above  satisfies the following  
\begin{align*}
  \sum_{j=1}^kn_ji_j& \leq \sum_{j=1}^kn_j(d-(k-j))=  n_1(d-(k-1))+\ldots + n_{k-1}(d-1)+n_kd\\
  &=d\sum_{j=1}^kn_j-\sum_{j=1}^{k-1}n_j(k-j) \leq d\n- (k-1)k/2
\end{align*}
  where in the second inequality we used $n_j\geq 1$ for all $j\leq k-1$. Hence we have equality if and only if 
  $i_j=d-(k-j)$ for $j\leq k$ and $n_j= 1$ for $j\leq k-1$.
  If these conditions are satisfied,  we furthermore have 
  $$
  \sum_{j=1}^kn_ji_j=dn-1 \quad \text{ if and only if }\quad \ k=2, 
  $$
  i.e. $i_1=d-1$, $i_2=d$,  $n_1=1$, $n_2=n-1$.
 Therefore the term
  $x^{d\n -1}$ appears in $g(x)^\n$ only once,    with coefficient equal to (a positive integer multiple of) $c_d^{\n-1}c_{d-1}$. Hence   $x^{\n d-1}$ appears
in  the left of \eqref{eq-dn} with coefficient $a_{0,\n} c_{d-1}c_d^{\n-1}$;  as this  coefficient must be zero and $a_{0,\n}c_d\neq 0$, we get  $c_{d-1}=0$. The induction base is proved.

Assume   $c_{d-1}=c_{d-2}=\ldots = c_{d-h+1}=0$.
We   have a non-zero coefficient of $x^{\sum_1^k n_ji_j}$
   only if $i_j\not\in\{  i_{d-1},i_{d-2},\ldots, i_{d-h+1}\}$, in which case we have
   \begin{align*}
  \sum_{j=1}^kn_ji_j& \leq n_1(d-h-k+2)+\ldots +n_{k-2}(d-h-1)+ n_{k-1}(d-h)+n_kd\\
  &=d\n- h \sum_{j=1}^{k-1}n_j- \sum_{j=1}^{k-2}j \leq d\n- h(k-1)-(k-1)(k-2)/2.
\end{align*}
The first inequality is an equality      if and only if $i_k=d$, and $i_{k-1}=d-h$, and $i_{j}=i_{j+1}-1$ for   $j<k-1$;
the second inequality is an equality      if and only if  $n_j=1$ for all $j\leq k-1$.
  If these conditions hold,   we furthermore have 
  $ 
  \sum_{j=1}^kn_ji_j=dn-h 
  $ if and only if $ k=2.$
  Therefore, arguing as above,   
  $x^{d\n -h}$ appears in $g(x)^\n$      with coefficient    $c_{d-h}c_d^{\n-1}$. Hence   it appears
  in the left    of \eqref{eq-dn} with coefficient $a_{0,\n} c_{d-h}c_d^{\n-1}$;  as this   must be zero, we get  $c_{d-h}=0$. The claim is proved.
  
  \
  
Thus  $C$ has equation $y=c_dx^d$,  hence $q$ is a $(1,d)$-point of $C$. 
But  $C$ is hypertangent to $B_3$ at $q$, hence Theorem~\ref{mirror-sm} implies that $B_3$ also has a $(1,d)$-point at $q$.

Assume $d\geq 3$.  
If $B$ is general  we can assume that no point in $B_3\cap (B_1\cup B_2)$ is a flex of $B_3$.  Hence $q$ is a $(1,2)$-point of $B_3$ and we get a contradiction. Therefore    $\Hyp_d(B,2)=\emptyset$ if $B$ is general.

If, instead, $B_3$ has a $(1,l)$-point in $q$  for some $l$ with $3\leq l\leq b_3$,   we get   $d=l$, hence $d$  is   determined by $B_3$, and $d\leq b_3$. Moreover   no term  of type $x^i$  with $i< d$ can appear in the equation of $B_3$ (for $q$ is a $(1,d)$-point), hence    $a_{i,0}=0$ for all $i<d$, and   
from \eqref{eq-dn} we derive   $$
c_dx^d+a_{d,0}x^d+  (\text{terms of degree }> d) =\lambda x^{d\n}.
$$
As $n\geq 2$ we get  $c_d=-a_{d,0}$. This proves that $C$ (if it exists) is uniquely determined by $B$ and $q$.

Assume $d=2$, then $q$ is a $(1,2)$-point for $C$ and   $B_3$. Now \eqref{eq-dn} gives
$$
c_2x^2+a_{2,0}x^2+a_{1,1} c_2x^{3}+ a_{3,0}x^3+ (\text{terms of degree }> 3) =\lambda x^{2\n}.
$$
We obtain 
      $a_{2,0}=-c_2$ and $a_{1,1}c_2=-a_{3,0}$, hence $a_{3,0}=a_{1,1}a_{2,0}$.
      Hence  
  $B_3$ is not a general curve of degree $\n$ (for its equation, \eqref{eq-B3}, must satisfy  $a_{3,0}=a_{1,1}a_{2,0}$).
Hence  $\Hyp_2(B,2)=\emptyset$ if $B$ is general. \eqref{mainc} is proved.
 
  If $B_3$ is not  general then, as before,    $C$ is     determined by the condition $c_2=-a_{2,0}$, hence it is   determined by $B$ and $q$. 
  
  Summarizing, for all $d\geq 2$ we proved that, for every  $q\in B_1\cap B_3$ there exists at most one curve $C\in \Hyp_d (B_2;p)  \cap \Hyp_d (B_3;q)$; as $q$ varies in $B_1\cap B_3$ we get at most $b_3$ curves  in $\Hyp(B,2)$.
 
  If $b_2=1$ the same argument applies by taking $q\in B_2\cap B_3$
  hence we might have $b_3$ new elements in $\Hyp (B,2)$.
  Hence  $|\Hyp_{\geq 2}(B,2)|\leq 2b_3$.
  
   If $b_2\geq 2$ we have $b_2$ choices for the point $p\in B_1\cap B_2$, hence   we obtain 
 $|\Hyp_{\geq 2}(B,2)|\leq b_2b_3$.
   \end{proof} 
For an example where the bound in the theorem is attained, see subsection~\ref{ss-sharp}.
\begin{cor}[of the proof]\label{cor:Vojta_bound}
If  $B$ is a 3C-curve  of degree at least $4$, and $C\in \Hyp(B,2)$, then $\deg C\leq b_3$
(hence $\deg C\leq \deg B-2$).
\end{cor}
It is worth stating the following  simple consequence of  our earlier results.
\begin{thm}\label{thm-3Cgen}
The   set $\Hyp(B,2)$  of a 
     3C-curve $B$ of degree at least $5$ is finite, and it is empty if $B$ is general. 
\end{thm} 
\begin{proof}
If $d=1$ this follows from 
 Proposition~\ref{prop-3C1}. Assume $d\geq 2$; then this
 follows from  
 Theorem~\ref{3C}
 if $b_1>1$, and  from
 Theorem~\ref{main} if $b_1=1$.\end{proof}

\section{Examples}
\label{sec-examples} 
We collect here some examples and special cases related to our   results.

\subsection{Examples of hypertangent curves}\label{ss-ex}
Let $B,C \subset \PP^2$ be two integral curves as in Theorem \ref{mirror-sm}, i.e. they have degree at least 2 and $B \cap C = \{ q \}$ with $q$ a unibranch point for $B$ and $C$. The following examples show that there are cases in which $q$ is a $(1,l)$-point for $B$ and $q$ is a singular point of $C$, i.e. is a $(m,lm)$ point with $m \geq 2$. In particular, in this setting, the type of the point $q$ for $B$ is different from its type as a point of $C$. 

(1) Consider the following two curves:
\[
  B :  y = x^2 \qquad C: (y - x^2)^3 + y^7.
\]
One easily  checks that $C\in \Hyp(B;q)$ where $q$ is the origin, and that $B$ has a $(1,2)$ point while $C$ has a $(3,6)$ point, so that $m=3$ and $l=2$ in our notation. More   generally, the following pair of curves,
\[
  B :  y = x^2 \qquad C: (y - x^2)^c + y^{2c+1},
\]
gives similar examples where $B$ has a $(1,2)$ point and $C$ has a $(c,2c)$ point for every $c \geq 3$.

(2) 
This second example  shows that a similar  situation can happen when $l$ is not  2. Consider the following two curves:
\[
  B: y = x^3 \qquad C: (y-x^3)^3 + y^9.
\]
In this case one checks that  $q$ is a $(1,3)$-point for $B$ while it is a $(3,9)$-point for $C$; as before  $C\in \Hyp(B;q)$.

We stress that both   these examples do not yield examples of curves in $\Hyp(\tilde{B},2)$ for a 3C-curve $\tilde{B}$ containing $B$. In fact, the proof of Theorem \ref{main} shows that the curves $C$ described above will have to meet the curve $\tilde B$ in more than two points.

\subsection{Explicit examples of hyper-bitangent conics}
\label{ss-sharp}
We show that the bounds in Theorem \ref{main} are sharp, in the sense that there are examples in which the set $\Hyp_{\geq 2}(B,2)$ consists of exactly $2b_3$ curves.

\begin{figure}[t!hht]  \centering 
 \includegraphics*[scale=.3]{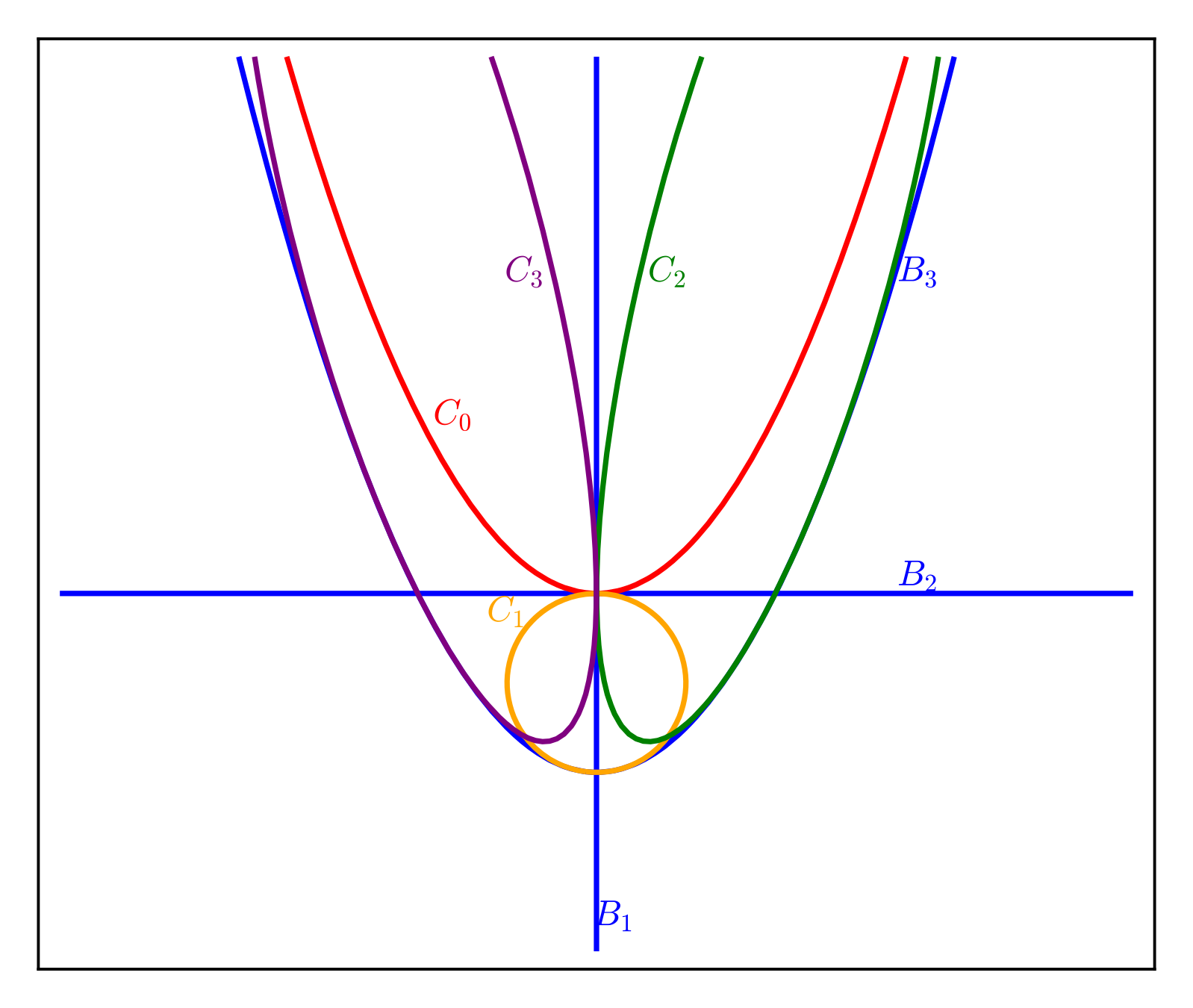}%
 \caption{Example of 4 conics in $\Hyp_2(B,2)$}\label{fig_4conic}
 \end{figure}

  For an explicit example consider the following three components of a 3C-curve $B = B_1 \cup B_2 \cup B_3$ with $(b_1,b_2,b_3) = (1,1,2)$:
  \[
    B_1: x = 0, \qquad B_2: z = 0, \qquad B_3 : zy = x^2 - y^2.
  \]
  Let $\{ q_1,q_2,q_3,q_4 \} = (B_1 \cap B_3) \cup (B_2 \cap B_3)$, where
  \[
    q_0 = (0:0:1), \quad q_1 = (0:-1:1), \quad q_2= (1:1:0),\quad q_3 = (1:-1:0).
  \]
  Then for $i=0,1,2,3$ the following conics $C_i$ satisfy $C_i \in \Hyp_2(B_3;q_i)$:
  \begin{align*}
    &C_0: zy = x^2, & C_1: zy = -z^2-x^2, \\
    &C_2: 4xz + 8xy - 8x^2 = z^2, & C_3: 4xz + 8xy + 8x^2 = -z^2.
  \end{align*}
  Moreover, $C_0, C_1 \in \Hyp_2(B_2;p)$ and $C_2,C_3 \in \Hyp_2(B_1;p)$ where $p = (0:1:0)$. In other words $B_2$ is the tangent line to $C_0$ and $C_1$ at $p$, while $B_1$ is the tangent line to $C_2$ and $C_3$ at $p$.  Therefore, for every $i = 0,1,2,3$, $C_i \in \Hyp_2(B;p,q_i)$ and hence $\lvert \Hyp_2(B,2) \rvert = 4$ reaching the upper bound of \eqref{maina} in the statement of Theorem \ref{main}.
 A picture of these conics in the affine patch $y=1$ can be seen in figure \ref{fig_4conic}.

\subsection{Curves with many components}
If a curve $B$ has more 
  than three components, it is not hard to prove that the only hyper-bitangent curves to $B$ are lines, and describe such lines 
 effectively. 
We  include the analysis of this case for completeness.
 \begin{prop}\label{cor-morecomp}
 Let $B$ be a reduced curve  of degree $b$  with $c\geq 4$  irreducible components, such that every point in the intersection of two components is a node of $B$. 
 
\begin{enumerate}[(a)]
\item
\label{cor-morecompa}
If $c\geq 5$ then  $\Hyp (B,2)=\emptyset$.
\item
 \label{cor-morecompb}
  If $d= 1$ then   $\Hyp_1(B,2)$ is finite; 
   moreover  
   \begin{enumerate} [(i)]
\item
 if $b=4$ then $|\Hyp_1(B,2)|=3$;
  \item
if $B$ does not contain two lines, then  $\Hyp_1(B,2) =\emptyset$;
   \item
if 
  $b\geq 5$ and $B$ is general then  $\Hyp_1(B,2)=\emptyset$.
\end{enumerate}
     \item
  \label{cor-morecompc}
 If  $d\geq 2$ then  $\Hyp_d(B,2)=\emptyset$.
 \end{enumerate}

\end{prop}
\begin{proof} Write $B=B_1\cup  \ldots \cup B_c$ with $c\geq 4$; by hypothesis through every node of $B$ there pass at most two components. 
Suppose there exists a curve  $C\in  \Hyp (B,2)$. Then $C$ meets $B$ in at most  two points, and   must intersect all  components of $B$. Hence $C\cap B=\{p,q\}$ and $p$, $q$ belong to exactly two components, say $p\in B_1\cap B_2$ and $q\in B_3\cap B_4$. In particular,  $B$ has only  4 components, proving \eqref{cor-morecompa}.

Let $d=1$  and   $C\in  \Hyp_1(B,2)$. As we said,  $C\cap B=\{p,q\}$. Since $B$ has finitely many (intersection) nodes,    $\Hyp_1(B,2)$ is finite.

If $b=4$ then $B$ is the union of 4 lines. It is clear that $\Hyp_1(B,2)$ is made of the three bitangent lines of $B$ not contained in $B$ (namely  the three lines not in $B$ and  joining a pair of nodes of $B$).

Let  $b\geq 5$.   We can assume that one component, $B_4$, has degree $\geq 2$.
Now, as $B_3$ and $B_4$ meet transversally in $q$ and $C\cap B_4=\{q\}$, the line   $C$ is  necessarily hypertangent to $B_4$ at $q$;
also, $B_3$ must have degree $1$, for it meets $C$ transversally in only one point.
Now, if such a line $C$  exists it is unique  and has to pass through the point $p$ as well. 
Arguing in the same way for $B_1$ and $B_2$ we have that at least one between $B_1$ and $B_2$ has degree $1$. Hence for  $ \Hyp_1(B,2)$ to be non-empty at least two components of $B$ have degree $1$.

Now, if the curve  $B$ is general, we can assume that no such   line exists, i.e.  we can assume that for every point
$q\in B_i\cap B_j$ the tangent lines  to $B_i$ and $B_j$ in $q$ do not  pass through any other intersection node of $B$.
This proves  that $\Hyp_1(B)$ is empty if $B$ is general. \eqref{cor-morecompb} is proved.

Let $d\geq 2$. By contradiction, let  $C\in \Hyp_d(B,2)$. As before, we   assume 
$C\cap B=\{p,q\}$ with $p\in B_1\cap B_2$ and $q\in B_3\cap B_4$.   Now,
  $B_1$ and $B_2$ meet transversally, hence $C$ in $p$ must be transverse to at least one between $B_1$ and $B_2$; say    $C$ is transverse to $B_1$.  Hence 
  $C$ must intersect $B_1$ in a further point, and this point must be $q$, which  is impossible as $q$ cannot belong to three components of $B$.
  \end{proof}

\subsection{Hypertangency  of rational curves}
\label{sec-snc} 
For an integral curve  $B\subset \PP^2$ of degree $b\geq 4$ having at most nodal singularities,
  Conjecture~\ref{conj:main} predicts that   there exist only  finitely many rational curves  hyper-bitangent to $B$, i.e. 
  the set $\cE(B)$ 
  is finite.
   We will provide an example showing  the necessity,  in the conjecture,    that $B$ have only nodal singularities. This example was already considered by Vojta in \cite[Example 3.5.3]{Vojta87}.

 For every integer $b\geq 3$  we  denote by  $Q_b\subset \PP^2$  the     curve given by the homogeneous equation 
 $$z^{b-1}y=x^b.$$  
 The curve $Q_b$ is smooth at $q_0=(0:0:1)$, with  tangent line, $L_{0}$, of equation $y=0$.
 We have 
 $Q_b\cap L_0=\{q_0\}$, 
so $L_0$ is   hypertangent to $Q_b$.
Next,  $Q_b$ has an $(b-1)$-fold unibranch point at $q_{\infty}=(0:1:0)$
with tangent line $L_{\infty}$ of equation $z=0$. We have 
 $Q_b\cap L_{\infty}=\{q_{\infty}\} $ 
so that  $L_{\infty}$ is hypertangent and $q_{\infty}$ is an $(b-1,b)$-point.
     \begin{center}
\begin{tikzpicture}[xscale=1.2,yscale=1.2]
 \draw[thick, blue](-0.2,3.3)to  [out=  0, in= 110]  (1.5,2.5);
  \draw[thick, blue](1.5,2.5)to  [out=  -60, in= 120]  (1,1.3); 
 \draw[thick, blue](0,2.1)to  [out= 0, in= 120]  (1,1.3);
 \draw[thick, blue](0,2.1)to  [out= 180, in= 0]  (-1,0.8); 
\draw[thick, blue](-1,0.8)to  [out=   180, in= 180]  (-0.2,3.3);

  \node at (-1.2,2.9){\color{blue}$ Q_b$};
   \draw [-] (-1.8,0.8)--(2,0.8)node at (-1.5,0.6){$L_0$};
   \draw [-] (0.25,2.8)--(1.4,0.5)node at (1.65,0.5){$L_{\infty}$};
  \node at (1.3,1.3){\color{blue}$q_{\infty}$};
  \filldraw[blue,  fill](0.99,1.3) circle (0.03);
  \node at (-1,1){\color{blue}$q_{0}$};
  \filldraw[blue,  fill](-1,0.8) circle (0.03);
 \end{tikzpicture}
  \end{center}
One checks easily that  $Q_b$ is integral,  has no other singular point, and is a rational curve; it  is, of course, not a nodal curve. 
Let us look at the set $\cE_d(Q_b)\subset \PP^{d(d+3)/2}$ of rational curves of degree $d$ hyper-bitangent to $Q_b$. If $d=1$ this has dimension 1, as it contains all lines through $q_{\infty}$. We have

\begin{prop}
For every $d\geq b$ and $b\geq 4$ we have
 \begin{equation}
 \label{eq-Qe1}
\dim \Hyp_d(Q_b,1)\geq 1 
\end{equation}
and 
\begin{equation}
 \label{eq-Qe2}
 \dim \cE_b(Q_b) \geq 1
\end{equation}
\end{prop}
 
\begin{proof}
Consider the curve $C_t$ of
 equation
$ 
y-x^b+ty^d=0  
$ 
with $t\in \C$.
 It is clear that $C_t$  lies in  $\Hyp_d(Q_b;q_0)$ for every $t\neq 0$   hence  \eqref{eq-Qe1} follows.
Notice that $C_t$ is smooth, hence not rational, for $t\neq 0$.

 To prove  \eqref{eq-Qe2} we will exhibit  a one-dimensional family of curves in $\cE_b(B)$.
 For every $t\neq 0,1$ let $R_t$ be the curve  having equation
$ 
z^{b-1}y=tx^b.
$ 
It is easy to check that $R_t$ is integral, and its singular locus consists of $q_{\infty}$ which is a     $(b-1,b)$-singular point, hence $R_t$ is  rational.
Moreover, we have
 $$
 R_t\cdot Q_b=bq_0+b(b-1)q_{\infty} 
 $$
 hence $R_t\in \Hyp (Q_b;q_0,q_{\infty})$. 
 As $t$ varies in $\C$ the curves $R_t$ form a family with a non-integral  member,  for $t=0$, hence the family is not constant. Therefore 
the $R_t$'s  give a one-dimensional subspace of $\cE_b(Q_b)$. \end{proof}

We stress that the condition that the integral curve $B \subset \PP^2$ has at most nodal singularities, although sufficient is not necessary. In other words there exist curves $B$, with worse singularities, for which $\cE(B)$ is still finite. A series of examples was constructed in \cite{CZ00}, see in particular the Proposition in page 2. To see the link with our setting, fixing $x,y,z$ as coordinates of $\PP^2$, the curve $B$ is given in op. cit. by the union of the line at infinity, the line $x = 0$, and a third component of affine equation $f(x,y) = 0$.

\bibliography{references}{}
\bibliographystyle{alpha}
\end{document}